\documentclass[12pt]{article}
\usepackage{latexsym,amsmath,amsthm,amssymb}
\usepackage{a4wide}
\usepackage{hyperref}
\usepackage{marginnote}
\usepackage{color}

\theoremstyle{plain}
\newtheorem{theorem}{Theorem}[section]

\newtheorem{lemma}[theorem]{Lemma}

\theoremstyle{definition}

\theoremstyle{remark}

\renewcommand{\thefootnote}{\arabic{footnote}}

\def\R{\mathbb R}

\def\va{\mathop{\rm Var}\nolimits}
\def\sign{\mathop{\rm sign}\nolimits}

\def\Om{\Omega}
\def\be{\beta}
\def\De{\Delta} 
\def\vphi{\varphi}
\def\na{\nabla}
\def\pa{\partial}
\def\la{\langle} 
\def\ra{\rangle} 
\def\lt{\left}
\def\rt{\right}
\def\o{\overline}

\numberwithin{equation}{section}

\title{A local proof of the dimensional Pr\'ekopa's theorem.}
\author{Van Hoang Nguyen\footnote{
School of Mathemacical Sciences, Tel Aviv University, Tel Aviv 69978, Israel. Email: vanhoang0610@yahoo.com}}

\begin{document}
\maketitle


\renewcommand{\thefootnote}{}

\footnote{Supported by a grant from the European Research Council}

\footnote{2010 \emph{Mathematics Subject Classification\text}: 26A51.}

\footnote{\emph{Key words and phrases\text}: Pr\'ekopa's theorem, Dimensional Pr\'ekopa's theorem, $L^2$-method}

\renewcommand{\thefootnote}{\arabic{footnote}}
\setcounter{footnote}{0}

\begin{abstract}
The aim of this paper is to find an expression for second derivative of the function $\phi(t)$ defined by
$$\phi(t) = \lt(\int_V \vphi(t,x)^{-\beta} dx\rt)^{-\frac1{\be -n}},\qquad \beta\not= n,$$
where $U\subset \R$ and $V\subset \R^n$ are open bounded subsets, and $\vphi: U\times V\to \R_+$ is a $C^2-$smooth function. As a consequence, this result gives us a direct proof of the dimensional Pr\'ekopa's theorem based on a local approach.
\end{abstract}

\section{Introduction}
The Pr\'ekopa's theorem \cite{Pre} says that marginals of log-concave functions are log-concave, i.e, if $\vphi:\R^{n+1}\to \R$ is convex, then the function $\phi$ defined by
\begin{equation}\label{eq:Pre}
\phi(t) = -\log\lt(\int_{\R^n} e^{-\vphi(t,x)} dx\rt)
\end{equation}
is convex on $\R$. By modifying $\vphi$ if necessary, we can replace $\R^{n+1}$ by any of its open convex subset $\Om$, and the integration in~\eqref{eq:Pre} is taken in the section $\Om(t) = \{x\in \R^n\, :\, (t,x)\in \Om\}$. The Pr\'ekopa's theorem is a direct consequence of the Pr\'ekopa-Leindler inequality which can be seen as the functional form of the Brunn-Minkowski inequality (see \cite{Gard}). The Brunn-Minkowski inequality is known to be one of the most important tools in analysis and geometry. It states that if $A, B$ are non-empty measurable subsets of $\R^n$ then
$$|A+ B|^{\frac1n}\geq |A|^{\frac1n} + |B|^{\frac1n},$$
where $A + B = \{a+b\, :\, a\in A,\, b\in B\}$ and $|\cdot|$ denotes the Lebesgue measure of the measurable set (See (\cite{BL,Gard,Led,Mau}) for the proofs and applications of the Brunn-Minkowski inequality).

A new proof of the Pr\'ekopa theorem is recently given in \cite{ABBN,BBN}. In these papers, the authors proved a local formulation for the second derivative of the function $\phi$ above. By using the convexity of $\vphi$, they show that $\phi''$ is nonnegative. This local approach was also used by D. Cordero-Erausquin (see \cite{Cordero}) to generalize a result of Berndtsson concerning the Pr\'ekopa's theorem for plurisubharmonic functions (see \cite{Berndtsson}). 

In this paper, we adapt the local approach given in \cite{ABBN,BBN} to find an expression for the second derivative of the function $\phi$ defined by
\begin{equation}\label{eq:phidimPrek}
\phi(t) = \lt(\int_V \vphi(t,x)^{-\beta} dx\rt)^{-\frac1{\be -n}},\qquad \beta\not= n,
\end{equation}
where $U\subset \R$ and $V\subset \R^n$ are open bounded subsets, boundary of $V$ is $C^\infty-$smooth, and $\vphi: U\times V\to \R_+$ is a $C^2-$smooth function on $U\times V$. For this purpose, we denote for each $t\in U$
$$d\mu_t = \frac{\vphi(t,x)^{-\be} dx}{\int_V \vphi(t,x)^{-\be} dx}$$
the probability measure on $V$. We also denote the corresponding symmetric diffusion operation with the invariant measure $\mu_t$ by
$$L_t u(x) = \De u(x) -\be \frac{\la\na_x \vphi(t,x),\na u(x)\ra}{\vphi(t,x)},$$
where $u$ is any function in $C^2(V)$. By using integration by parts, we have
$$\int_V L_tu(x) v(x) d\mu_t(x) = -\int_V\la \na u(x),\na v(x)\ra d\mu_t(x) + \int_{\pa V} v(x) \frac{\pa u}{\pa \nu}(x) d\mu_t(x),$$
where $\nu(x)=(\nu_1(x), \cdots, \nu_n(x))$ is the outer normal vector to $x\in \pa V$.

Since $\pa V$ is $C^\infty-$smooth, then $\nu$ is $C^\infty-$smooth on $\pa V$ and it can be extended to a $C^\infty-$smooth map on a neighbourhood of $\pa V$. Hence the second fundamental form $II$ of $\pa V$ at $x\in \pa V$ is defined by
$$II_x(X,Y) = \sum_{i,j=1}^n X_iY_j\pa_i(\nu_j)(x),$$
for any two vector fields $X =(X_1,\cdots,X_n)$ and $Y =(Y_1,\dots, Y_n)$ in $\pa V$.

In the sequel, we denote by $\nabla f$ and $\nabla^2f$ the gradient and Hessian matrix of a function $f$, respectively. We also denote by $\|\cdot\|_{HS}$ the Hilbert-Schmidt norm on the space of square matrices. When $f$ is function of the variables $t$ and $x$, we write $\nabla_x f$ and $\nabla^2_x f$ for the gradient and Hessian matrix of $f$ which are taken only on $x$, respectively.

Our first main theorem of this paper is the following:

\begin{theorem}\label{maintheorem}
Suppose that $V$ has $C^\infty-$smooth boundary, and $\vphi$ is $C^\infty-$smooth up to boundary of $U\times V$. Let $\phi$ be defined by~\eqref{eq:phidimPrek} then
\begin{align}\label{eq:secondderivative}
\frac{\phi''(t)}{\phi(t)} &= \frac{\be}{\be -n}\int_V \frac{\la(\na^2_{(t,x)}\vphi )X,X\ra}\vphi d\mu_t +\frac{\be^2}{\be -n}\int_V\lt(\|\na^2 u\|_{HS}^2 -\frac1n (\De u)^2\rt) d\mu_t\notag\\
&\quad + \frac{\be}{|\be -n|}\int_V \lt(\sqrt{\frac{|\be-n|}n}\De u -\mathrm{sign}(\be-n)\sqrt{\frac{n}{|\be -n|}} \int_V \frac{\pa_t\vphi}{\vphi} d\mu_t\rt)^2 d\mu_t \notag\\
&\quad + \frac{\be^2}{\be -n} \int_{\pa V} II(\na u,\na u) d\mu_t,
\end{align}
where $u$ is the solution of the equation
\begin{equation}\label{eq:Eeq}
L_t u = \frac{\pa_t \vphi(t,\cdot)}{\vphi(t,\cdot)} - \int_V \frac{\pa_t\vphi(t,x)}{\vphi(t,x)} d\mu_t(x)\quad\text{and}\quad \frac{\pa u(x)}{\pa\nu(x)} =0,\quad x\in \pa V,
\end{equation}
and $X$ denotes the vector field $(1,\be\na u(x))$ in $\R^{n+1}$. 
\end{theorem}
Since $\frac{\pa u(x)}{\pa \nu(x)} =0$ for every $x\in \pa V$, hence $\nabla u(x)\in T_x(\pa V)$ (the tangent space to $\pa V$  at $x\in \pa V$). This implies that $II(\nabla u,\nabla u)$ is well-defined on $\pa V$. Theorem \ref{maintheorem} is proved in the next section. We will need the following classical fact about the existence of the solution of the elliptic partial differential equation (see \cite{KM} and references therein): 
\begin{lemma}\label{solution}
If $V$ has $C^\infty-$smooth boundary $\pa V$, and $\vphi$ is $C^\infty$-smooth up to boundary of $V$, then for any function $f\in C^\infty(\overline{V})$, $\int_V f(x) d\mu_t(x) = 0$ there exists a function $u\in C^{\infty}(\overline{V})$ such that $L_tu = f$ and $\frac{\pa u(x)}{\pa\nu(x)} = 0$ on $\pa V$. 
\end{lemma}

Our second main theorem of this paper is the dimensional Pr\'ekopa's theorem which is considered as a direct consequence of Theorem~\ref{maintheorem} and stated in the following theorem. The first part of this theorem concerns the convex case, and the second part concerns the concave case.

\begin{theorem}\label{DimPrek}
Let $\Om\subset \R^{n+1}$ be a convex open subset, and let $\vphi: \Om\to \R_+$ be a $C^2-$smooth function up to boundary of $\Om$. For $t\in \R$, we define the section $\Om(t) =\{x\in \R^n\, :\, (t,x)\in \Om\}$. Then the following assertions hold:
\begin{description}
\item (i) If $\vphi$ is convex on $\Om$, and $\be >n$, then the function $\phi$ defined by
$$\phi(t) = \lt(\int_{\Om(t)} \vphi(t,x)^{-\be} dx\rt)^{-\frac1{\beta -n}},$$
is convex on $\R$.
\item (ii) If $\vphi$ is concave on $\Om$, and $\be > 0$, then the function $\phi$ defined by
$$\phi(t) = \lt(\int_{\Om(t)} \vphi(t,x)^{\be} dx\rt)^{\frac1{\beta +n}},$$
is concave on $\R$.
\end{description}
\end{theorem}
Finally, we remark that the Pr\'ekopa's theorem can be deduced from Theorem~\ref{DimPrek} by letting $\be$ tend to infinity since
$$\lim\limits_{\be\to\infty}\,(\be -n)\lt[\lt(\int_{\Om(t)} \lt(1+ \frac{\vphi(t,x)}{\be}\rt)_+^{-\be} dx\rt)^{-\frac1{\be-n}} -1\rt]= -\log\lt(\int_{\Om(t)} e^{-\vphi(t,x)} dx\rt)$$
and
$$\lim\limits_{\be\to\infty}\,(\be +n)\lt[\lt(\int_{\Om(t)} \lt(1- \frac{\vphi(t,x)}{\be}\rt)_+^{\be} dx\rt)^{\frac1{\be+n}} -1\rt]= \log\lt(\int_{\Om(t)} e^{-\vphi(t,x)} dx\rt),$$
where $a_+ = \max\{a,0\}$ denotes the positive part of $a$.


\section{Proof of main theorems}

We begin this section by giving the proof of Theorem~\ref{maintheorem}. Our proof is direct and similar the method used in \cite{NVH}.

\emph{\bf Proof of Theorem~\ref{maintheorem}\text}: If $\be = 0$, then~\eqref{eq:secondderivative} is evident since $\phi$ is a constant function.

If $\be\not=0$, then~\eqref{eq:secondderivative} is equivalent to
\begin{align}\label{eq:secondderivative1}
\frac{\beta -n}{\beta}\frac{\phi''(t)}{\phi(t)} &= \int_V \frac{\la(\na^2_{(t,x)}\vphi) X,X\ra}\vphi d\mu_t +\be\int_V\lt(\|\na^2 u\|_{HS}^2 -\frac1n (\De u)^2\rt) d\mu_t\notag\\
&\quad + \sign\lt(\be -n\rt)\int_V \lt(\sqrt{\frac{|\be-n|}n}\De u -\sign(\be-n)\sqrt{\frac{n}{|\be -n|}} \int_V \frac{\pa_t\vphi}{\vphi} d\mu_t\rt)^2 d\mu_t \notag\\
&\quad + \be \int_{\pa V} II(\na u,\na u) d\mu_t.
\end{align}
By a direct computation, we easily get 
\begin{align}\label{eq:daohamcaphai}
\frac{\beta -n}{\beta}\frac{\phi''(t)}{\phi(t)} &= \int_V \frac{\pa^2_{tt}\vphi(t,x)}{\vphi(t,x)} d\mu_t(x) -(\beta +1)\va_{\mu_t} \lt(\frac{\pa_t\vphi(t,\cdot)}{\vphi(t,\cdot)}\rt) \notag\\
&\qquad + \frac{n}{\be-n}\lt(\int_V \frac{\pa_t\vphi(t,x)}{\vphi(t,x)} d\mu_t(x)\rt)^2,
\end{align}
where $\va_{\mu_t}(f) := \int_V f^2 d\mu_t - (\int_V f d\mu_t)^2$ denotes the variance of any function $f$ on $V$ with respect to $\mu_t$.

Let $u\in C^\infty(\overline{V})$ be the solution of the equation~\eqref{eq:Eeq}.
Since $\mu_t$ is a probability measure on $V$, then we have
$$\va_{\mu_t}\lt(\frac{\pa_t\vphi(t,\cdot)}{\vphi(t,\cdot)}\rt) = -\int_V (L_t u)^2 d\mu_t + 2\int_V \lt(\frac{\pa_t\vphi}{\vphi}-\int_V \frac{\pa_t \vphi}\vphi d\mu_t\rt) L_t u\, d\mu_t.$$
 Using integration by parts and the fact $\int_V L_t u \,d\mu_t = 0$, we get
\begin{equation}\label{eq:integrationbyparts}
\int_V \lt(\frac{\pa_t\vphi}{\vphi} -\int_V \frac{\pa_t \vphi}\vphi d\mu_t\rt) L_t u\, d\mu_t = -\int_V \frac{\la \na_x(\pa_t\vphi),\na u\ra}{\vphi} d\mu_t + \int_V \frac{\pa_t\vphi}{\vphi}\frac{\la\na_x\vphi,\na u\ra}{\vphi}d\mu_t.
\end{equation}
It follows from integration by parts (see also the proof of Theorem $1$ in \cite{NVH}) that
\begin{align}\label{eq:intLu^2}
\int_V (L_t u)^2 d\mu_t & = \int_V \|\na^2 u\|_{HS}^2 d\mu_t + \be\int_V \frac{\la(\na_x^2\vphi)\na u,\na u\ra}\vphi d\mu_t \notag\\
&\qquad - \be \int_V \frac{\la \na_x\vphi,\na u\ra^2}{\vphi^2} d\mu_t -\int_{\pa V}\la (\na^2 u)\na u,\nu\ra d\mu_t. 
\end{align}

From~\eqref{eq:integrationbyparts} and~\eqref{eq:intLu^2}, we get an expression of $\va_{\mu_t}(\pa_t\vphi/\vphi)$ as follows
\begin{align}\label{eq:variance}
\va_{\mu_t}\lt(\frac{\pa_t\vphi(t,\cdot)}{\vphi(t,\cdot)}\rt)& =-2\int_V \frac{\la \na_x(\pa_t\vphi),\na u\ra}{\vphi} d\mu_t + 2\int_V \frac{\pa_t\vphi}{\vphi}\frac{\la\na_x\vphi,\na u\ra}{\vphi}d\mu_t \notag\\
&\qquad -\int_V \|\na^2 u\|_{HS}^2 d\mu_t - \be\int_V \frac{\la(\na_x^2\vphi)\na u,\na u\ra}\vphi d\mu_t \notag\\
&\qquad + \be \int_V \frac{\la \na_x\vphi,\na u\ra^2}{\vphi^2} d\mu_t +\int_{\pa V}\la (\na^2 u)\na u,\nu\ra d\mu_t.
\end{align}
It follows from the definition of $L_t$ that
\begin{equation}\label{eq:dem}
\be^2 \int_V \frac{\la \na_x\vphi,\na u\ra^2}{\vphi^2} d\mu_t = \int_V\lt[(L_t u)^2 +(\De u)^2\rt]d\mu_t -2\int_V \De u\, L_t u\, d\mu_t.
\end{equation}
Plugging~\eqref{eq:intLu^2} and~\eqref{eq:Eeq} into~\eqref{eq:dem}, we obtain
\begin{align}\label{eq:carer}
\be(\be+1)\int_V \frac{\la \na_x\vphi,\na u\ra^2}{\vphi^2} d\mu_t &= \int_V \|\na^2 u\|_{HS}^2 d\mu_t + \int_V (\De u)^2 d\mu_t\notag\\
&\quad +\be\, \int_V \frac{\la(\na_x^2\vphi)\na u,\na u\ra}\vphi d\mu_t - 2\int_V \frac{\pa_t\vphi}{\vphi}\,\De u\, d\mu_t\notag\\
& \quad + 2\int_V\De u\lt(\int_V \frac{\pa_t \vphi}\vphi d\mu_t\rt) d\mu_t -\int_{\pa V} \la (\na^2 u)\na u,\nu\ra d\mu_t.
\end{align}
Moreover, using again integration by parts, we have
\begin{align*}
\int_V \frac{\pa_t\vphi}{\vphi}\frac{\la\na_x\vphi,\na u\ra}{\vphi}d\mu_t &= -\frac1\be\int_V \frac{\pa_t\vphi(t,x)}{\vphi(t,x)}\la\na_x(\vphi(t,x)^{-\be}),\na u(x)\ra \, dx\\
&=\frac1\be\int_V\frac{\la\na_x(\pa_t\vphi),\na u\ra}\vphi d\mu_t -\frac1\be\int_V\frac{\pa_t\vphi}{\vphi}\frac{\la\na_x\vphi,\na u\ra}{\vphi}d\mu_t\\
&\quad + \frac1\be\int_V \frac{\pa_t\vphi}\vphi \,\De u \, d\mu_t.  
\end{align*}
Or, equivalent
\begin{equation}\label{eq:IBP}
(\be +1)\int_V \frac{\pa_t\vphi}{\vphi}\frac{\la\na_x\vphi,\na u\ra}{\vphi}d\mu_t = \int_V\frac{\la\na_x(\pa_t\vphi),\na u\ra}\vphi d\mu_t + \int_V \frac{\pa_t\vphi}\vphi \,\De u \, d\mu_t.
\end{equation}
Plugging~\eqref{eq:variance},~\eqref{eq:carer}, and~\eqref{eq:IBP} into~\eqref{eq:daohamcaphai}, we obtain 
\begin{align}\label{eq:remainboundary}
\frac{\be -n}\beta\frac{\phi''(t)}{\phi(t)} &=
\int_V \frac{\pa^2_{tt}\vphi}{\vphi} d\mu_t +2\be \int_V \frac{\la \na_x(\pa_t\vphi),\na u\ra}{\vphi} d\mu_t + \be^2\int_V \frac{\la(\na_x^2\vphi)\na u,\na u\ra}\vphi d\mu_t\notag\\
&\quad +\be \int_V \|\na^2u\|_{HS}^2 d\mu_t -2\int_V \De u\lt(\int_V \frac{\pa_t \vphi}\vphi d\mu_t\rt)\, d\mu_t -\int_V (\De u)^2 d\mu_t \notag\\
&\quad + \frac{n}{\be -n} \lt(\int_V \frac{\pa_t\vphi(t,x)}{\vphi(t,x)} d\mu_t(x)\rt)^2 -\be \int_{\pa V}\la (\na^2 u)\na u,\nu\ra d\mu_t\notag\\
&=\int_V \frac{\pa^2_{tt}\vphi}{\vphi} d\mu_t +2\be \int_V \frac{\la \na_x(\pa_t\vphi),\na u\ra}{\vphi} d\mu_t + \be^2\int_V \frac{\la(\na_x^2\vphi)\na u,\na u\ra}\vphi d\mu_t\notag\\
&\quad +\be\int_V\lt(\|\na^2 u\|_{HS}^2 -\frac1n (\De u)^2\rt) d\mu_t+\frac{\be -n}n \int_V (\De u)^2 d\mu_t\notag\\
&\quad -2\int_V \De u\lt(\int_V \frac{\pa_t \vphi}\vphi d\mu_t\rt)\, d\mu_t +\frac{n}{\be -n} \lt(\int_V \frac{\pa_t\vphi(t,x)}{\vphi(t,x)} d\mu_t(x)\rt)^2\notag\\
&\quad -\be \int_{\pa V}\la (\na^2 u)\na u,\nu\ra d\mu_t.
\end{align}
To finish our proof, we need to treat the term on boundary in~\eqref{eq:remainboundary}. Since $\frac{\pa u}{\pa \nu} = 0$ on $\pa V$, then $\na u(x) \in T_x(\pa V)$ for every $x\in \pa V$, and
\begin{equation}\label{eq:2ndform}
\la (\na^2u(x))\na u(x),\nu(x)\ra = -II_x(\na u(x),\na u(x)),\qquad x\in \pa V.
\end{equation}
Combining~\eqref{eq:remainboundary} and~\eqref{eq:2ndform}, and denoting $X(t,x) = (1,\be\na u(x))$ with $(t,x)\in U\times V$, we get~\eqref{eq:secondderivative1}. Then Theorem~\ref{maintheorem} is completely proved.

In the following, we use Theorem~\ref{maintheorem} to prove the dimensional Pr\'ekopa's theorem (Theorem~\ref{DimPrek}).

\emph{\bf Proof of Theorem~\ref{DimPrek}\text}: By using an approximation argument, we can assume that $\Om$ is bounded and $\vphi$ is $C^\infty$-smooth up to boundary of $\Om$.

\emph{Part $(i)$\text}: We first prove when $\Om = U\times V$ with $U\subset \R$, and $V\subset \R^n$ has $C^\infty$-smooth boundary $\pa V$. Since $\be > n$, then applying Theorem~\ref{maintheorem}, we have
\begin{align}\label{eq:pos}
\frac{\be -n}\be\frac{\phi''(t)}{\phi(t)} &=\int_V \frac{\la(\na^2_{(t,x)}\vphi) X,X\ra}\vphi d\mu_t +\be \int_V\lt(\|\na^2 u\|_{HS}^2 -\frac1n (\De u)^2\rt) d\mu_t\notag\\
&\quad + \int_V \lt(\sqrt{\frac{\be-n}n}\De u -\sqrt{\frac{n}{\be -n}} \int_V \frac{\pa_t\vphi}{\vphi} d\mu_t\rt)^2 d\mu_t\notag\\
&\quad + \be  \int_{\pa V} II(\na u,\na u) d\mu_t,
\end{align}
where $II$ denotes the second fundamental form of $\pa V$, and $u$ is the $C^\infty-$smooth solution of the equation~\eqref{eq:Eeq} with $L_t = \De -\be \la \na_x\vphi,\cdot\ra/\vphi$, and $X$ denotes the vector field $(1,\be\na u)$ in $\R^{n+1}$.

We have $II_x(\na u(x),\na u(x)) \geq 0$, $x\in \pa V$ because of the convexity of $V$. By Cauchy-Schwartz inequality, we have 
$$\|\na^2 u\|_{HS}^2 \geq \frac1n\, (\De u)^2.$$
As a consequence of the convexity of $\vphi$, we obtain $\na^2 \vphi \geq 0$ in the sense of symmetric matrix. All the integrations on the right hand side of~\eqref{eq:pos} hence are nonnegative. This implies that $\phi'' \geq 0$, or $\phi$ is convex.

In the general case, there exists an increasing sequences of $C^\infty$-smooth open convex $\Om_k$ such that
$$\Om_k = \{(t,x)\, :\, \rho_k(t,x) < 0\},$$
with $\rho_k \in C^\infty(\R^{n+1})$, $k= 1,2\cdots$ are convex functions, and $\Om =\bigcup_k \Om_k$. Hence, by using an approximation argument, we can assume that
$$\Om =\{(t,x)\, :\, \rho(t,x) < 0\}$$
with a $C^{\infty}-$smooth convex function $\rho$, and $\vphi$ is defined in a neighborhood of $\Om$. Since the convexity is local, it is enough to prove that $\phi$ is convex in a neighborhood of each $t$. Fix $t_0$, choose a small enough neighborhood $U$ of $t_0$ such that
$$(U\times \R^n) \cap \o{\Om} \subset U\times V$$
and $\rho$, $\vphi$ are defined in $U\times V$, where $V$ is convex subset of $\R^n$ and has $C^\infty-$smooth boundary $\pa V$. Define $\rho_0 = \max\{\rho,0\}$, then $\rho_0$ is a convex function in $U\times V$. With $N >0$, we know that the function
$$\phi_N(t) = \lt(\int_V \lt(\vphi(t,x) + N\rho_0(t,x)\rt)^{-\be} dx\rt)^{-\frac1{\beta-n}}$$
is convex in $U$. Moreover, $\phi_N(t)\to \phi(t)$ in $U$ as $N$ tends to infinity, then $\phi$ is convex in $U$. This finishes the proof of the convexity of $\phi$.

\emph{Part $(ii)$\text}: As explained in the proof of the part $(i)$ above, it suffices to prove the part $(ii)$ in the case $\Om =U\times V$ with $U\subset\R$ and $V\subset \R^n$ are bounded open convex subsets, and $\pa V$ is $C^\infty-$smooth. Since $\be > 0$, by applying Theorem~\ref{maintheorem} to $-\beta$ instead of $\beta$, we have
\begin{align}\label{eq:neg}
\frac{\be +n}\be\frac{\phi''(t)}{\phi(t)} &=\int_V \frac{\la(\na^2_{(t,x)}\vphi) X,X\ra}\vphi d\mu_t -\be \int_V\lt(\|\na^2 u\|_{HS}^2 -\frac1n (\De u)^2\rt) d\mu_t\notag\\
&\quad - \int_V \lt(\sqrt{\frac{\be+n}n}\De u +\sqrt{\frac{n}{\be +n}} \int_V \frac{\pa_t\vphi}{\vphi} d\mu_t\rt)^2 d\mu_t\notag\\
&\quad - \be \int_{\pa V} II(\na u,\na u) d\mu_t.
\end{align}
Using the arguments in the proof of part $(i)$ and the concavity of $\vphi$, we get $\phi''(t)\leq 0$ from~\eqref{eq:neg}, or $\phi$ is concave.

\subsection*{Acknowledgment}
The author would like to sincerely thank anonymous referee for many useful and valuable comments which improved the quality of this paper.

\end{document}